\documentclass[12pt]{article}
\usepackage{amssymb, url}
\usepackage{amsmath,amsthm, xpatch}
\usepackage{mathtools}
\usepackage{dsfont}
\usepackage{verbatim}
\usepackage{graphicx}
\usepackage{epstopdf}
\usepackage{fullpage}
\usepackage[hang,flushmargin]{footmisc}
\usepackage{tikz}
\usepackage{enumitem, hyperref}
\usepackage{datetime}
\usepackage{bbm}
\usepackage{color}
\usetikzlibrary{calc}
\usetikzlibrary{decorations.pathreplacing,calligraphy}

\newtheorem{theorem}{Theorem}[]

\newtheorem{lemma}{Lemma}
\newtheorem{conjecture}{Conjecture}

\newtheorem{claim}{Claim}[theorem]

\makeatletter
\newcommand{\newreptheorem}[2]{\newtheorem*{rep@#1}{\rep@title}\newenvironment{rep#1}[1]{\def\rep@title{#2 \ref*{##1}}\begin{rep@#1}}{\end{rep@#1}}}
\makeatother
\newreptheorem{theorem}{Theorem}

\theoremstyle{plain}

\newcounter{caseCount}
\AtBeginEnvironment{claim}{\setcounter{caseCount}{0}}

\providecommand{\keywords}[1]
{
  \small	
  \textbf{\textit{Keywords---}} #1
}

\begin{document}
\title{Connected equitably $\Delta$-colorable realizations with $k$-factors.}
\author{James M.\ Shook$^{1,2}$}

\footnotetext[1]{National Institute of Standards and Technology, Computer Security Division, Gaithersburg, MD; {\tt james.shook@nist.gov}.}
\footnotetext[2]{This article is a U.S. Government work and is not subject to copyright in the USA.}
\maketitle
\begin{abstract}
A graph $G$ is said to be equitably $c$-colorable if its vertices can be partitioned into $c$ independent sets that pairwise differ in size by at most one. Chen, Lih, and Wu conjectured that every connected graph $G$ with maximum degree $\Delta(G)\geq 2$ has an equitable coloring with $\Delta(G)$ colors, except when $G$ is complete, an odd cycle, or a balanced bipartite graph with odd sized partitions. Suppose $G$ is a connected graph with a $k$-factor (a regular spanning subgraph) $F$ such that $G$ is not complete, a $1$-factor, nor an odd cycle. When $k\geq 1$ we demonstrate that there is a connected $(k-1)$ edge-connected equitably $\Delta(G)$-colorable graph $H$ with a $k$-factor $F'$ such that $G-E(F)=H-E(F')$. If we drop the requirement that $G-E(F)=H-E(F')$, then we can say more. Considering the non-increasing degree sequence $\pi=(d_{1},\ldots, d_{n})$ of $G$ where $d_{i}=deg_{G}(v_{i})$ for all vertices $\{v_{1},\ldots,v_{n}\}$ of $G$, we call $m(\pi)=\max\{i|d_{i}\geq i\}$ the strong index of $\pi$. For $k\geq 0$, we can show that for every $$c\geq \max_{l\leq m(\pi)}\bigg\{\bigg\lfloor\frac{d_{l}+l}{2}\bigg\rfloor\bigg\}+1$$ we can find a connected $(k-1)$ edge-connected equitably $c$-colorable realization $H$ of $\pi$ that has a $k$-factor. In a third theorem we show that if  $d_{d_{1}-d_{n}+1}\geq d_{1}-d_{n}+k-1$, then some realization of $\pi$ has a $k$-factor. Together, these three theorems allow us to prove that for all $k$, there is a connected equitably $\Delta(G)$-colorable realization $H$ of $\pi$ with a $k$-factor. Thus,  giving support to the validity of the Chen-Lih-Wu Conjecture.
\end{abstract}


\section{Introduction}
All graphs in this paper are finite and simple. The $n$ vertex graphs $C^{n}$, $I^{n}$, $K^{n}$, and $K^{a,b}$ denote a cycle, an independent set, a complete graph, and a complete bipartite graph with $a+b=n$, respectively. We let $H*G$ denote the join of $H$ with $G$.

Let $G=(V,E)$ be a graph with $X,Y\subseteq V$. The induced graph on $X$ is denoted by $G[X]$, and we let $\overline{X}=V-X$. We denote $E_{G}(X,Y)$ to be the set of all edges of $G$ that have one end in $X$ and the other in $Y$, and let $e_{G}(X,Y)=|E_{G}(X,Y)|$ and $e_{G}(x,Y)=e_{G}(\{x\},Y)$.  The set of all vertices in $X$ that are adjacent in $G$ to vertices in $\overline{X}$ is represented by $\Gamma_{G}(X)$.  The edge-connectivity of $G$ is denoted by $\lambda(G)$, and $deg_{G}(v)$ denotes the degree of $v\in V$. The maximum and minimal degree of $G$ is denoted by $\Delta(G)$ and $\delta(G)$, respectively. If $\delta(G)=\lambda(G)$, then $G$ is said to be maximally edge-connected. A vertex coloring of $G$ using $c$ colors is a function $f: V\rightarrow \{1,\ldots,c\}$. We will assume all vertex colorings are proper in the sense that $f(x)\neq f(y)$ for all edges $xy\in E(G)$. We say a graph is $c$-colorable if it can be colored with $c$ colors. The color classes of $f$ is a partition $\{Y_{1}(f),\ldots, Y_{c}(f)\}$ of $V$ such that $f(v)=i$ for every $i\in \{1,\ldots, c\}$ and $v\in Y_{i}(f)$. An \textit{equitable $c$-coloring} of $G$ is a $c$-coloring $f$ in which  $||Y_{i}(f)|-|Y_{j}(f)||\leq 1$ for all $\{i,j\}\subseteq \{1,\ldots,c\}$. In this paper, we will assume that $|Y_{1}(f)|\geq \ldots \geq |Y_{c}(f)|$. See \cite{Lih2013} for a survey on equitable coloring. 


Given a graph $G=(V,E)$ with $V=\{v_{1},\ldots,v_{n}\}$ we call $(deg_{G}(v_{1}),\ldots, deg_{G}(v_{n}))$ a degree sequence of $G$. If a sequence of $n$ non-negative integers $\alpha=(d_{1},\ldots, d_{n})$ is a degree sequence of some graph, then we say $\alpha$ is graphic. All such graphs are said to realize $\alpha$, and we denote the set of all realizations of $\alpha$ by $\mathcal{R}(\alpha)$. Note that if $d'_{i}=d_{\sigma(i)}$ for some bijection $\sigma: [n]\rightarrow [n]$, then $\mathcal{R}(\alpha)=\mathcal{R}((d'_{1},\dots,d'_{n}))$. Thus, we conveniently assume all graphic sequences in this paper are non-increasing. If every $d_{i}$ of $\alpha$ is non-zero, then we say $\alpha$ is positive. We denote the degree sequence of a graph $G$ by $\pi(G)$, and  we let $\mathcal{R}(\pi(G),H)\subseteq \mathcal{R}(\pi(G))$ be the set of all $W\in \mathcal{R}(\pi(G))$ such that $E(H)\subseteq E(G)\cap E(W)$.

The rest of the paper is organized as follows. In Section~\ref{sec:motivation} we present our motivations for this work, and our main results are presented in Section \ref{sec:mainResults}. In Section \ref{sec:PastWork}, we review our previous work that lays the groundwork for this paper. In particular, Section \ref{sec:edgeconnectivity} discusses modifications to our theorem on edge connectivity, as presented in \cite{Shook2024}, to accommodate equitable coloring. The proofs of our main results can be found in the subsequent sections.

\subsection{Motivation}\label{sec:motivation}

In \cite{Kierstead2008}, Kierstead and Kostochka presented a short proof of a conjecture on equitable coloring that was first posed by Erd\H{o}s in 1964 and originally proved by Hajnal and Szemer\'{e}di. 

\begin{theorem}[Hajnal-Szemer\'{e}di \cite{Hajnal1970}]\label{thm:HSTheorem} For every positive integer r, each graph $G$ with $\Delta(G)\leq r$ has an equitable $(r+1)$-coloring.
\end{theorem}

This theorem cannot be improved in general since complete graphs and graphs with an odd cycle cannot be colored with $\Delta(G)$ colors. However, Brooks proved these are the only obstacles.

\begin{theorem}[Brooks's Theorem \cite{Brooks1941}]Every graph $G$ can be colored with $\Delta(G)$ colors unless $G$ contains a $K^{\Delta(G)+1}$ or $\Delta(G)=2$ and $G$ contains an odd cycle.
\end{theorem}

Later, Chen, Lih, and Wu conjectured that Brooks's theorem could be generalized to equitable $\Delta(G)$ colorings for connected graphs.

\begin{conjecture}[The Chen-Lih-Wu Conjecture \cite{Chen1994}]\label{con:Equitable}Every connected graph $G$ with maximum degree $\Delta(G)\geq 2$ has an equitable coloring with $\Delta(G)$ colors, except when $G$ is $K^{m}$, $C^{2m+1}$, or $G=K^{2m+1,2m+1}$ for all $m\geq 1$.
\end{conjecture}

The Chen-Lih-Wu conjecture is still open, but Kierstead and Kostochka proved it is true when $\Delta(G)\leq 4$ in \cite{Kierstead2012} or when $\Delta(G)\geq n/4$ in \cite{Kierstead2015}. The conjecture was proved true for bipartite graphs in \cite{Lih1996}, for split graphs in \cite{Chen1996}, and for many other special classes. Moreover, Chen, Lih, and Wu in \cite{Chen1994} showed that the conjecture is valid if it is valid for regular graphs.

Ultimately, as a consequence of our main results, we will show there is a realization of $\pi(G)$ that satisfies the Chen-Lih-Wu conjecture. However, for a graph $G$ and some $c\leq \Delta(G)$, we want to know if there is a realization of $\pi(G)$ that is not only equitably $c$-colorable but preserves select features of $G$. In particular, we focus on maintaining edge-connectivity and the existence of $k$-factors ($k$-regular spanning subgraphs).

\subsection{Main Results}\label{sec:mainResults}

Our first main result shows that if a graph $G$ has a $k$-factor and is not a forbidden graph, then we can find a realization that is not only equitably $\Delta(G)$-colorable, but has a $k$-factor, is connected, and contains a large spanning subgraph of $G$.

\begin{theorem}\label{thm:MaxDegreeEquitable}For $k\geq 1$, if $G=(V,E)$ is a simple $n$ vertex graph such that $G$ has a $k$-factor $F$ and $G$ is not a complete graph, a $1$-factor, nor a $2$-factor when $n$ is odd, then there is a realization $H\in \mathcal{R}(\pi(G), G-E(F))$ that is equitably $\Delta(G)$-colorable. In addition, if $G$ has at least $n-1$ edges, then we can require $H$ to be $(\delta(G)-1)$ edge-connected or $\delta(G)$ edge-connected when $\delta(G)$ is even or one. 
    
\begin{proof}See Section~\ref{sec:chenlihWu}.
\end{proof}
\end{theorem}

If we don't require $G-E(F)$ to be a subgraph of $H$, then we can show that there are realizations of $G$ that are equitably $\gamma$-colorable for some $\gamma\leq \Delta(G)+1$. We denote by $m(\pi)$ the largest index of $\pi(G)$ such that $d_{m(\pi)}\geq m(\pi)$. This number is sometimes called the strong index \cite{Zverovich1992} or the Durfee number \cite{Barrus2022}.

\begin{theorem}\label{thm:ExistDegreeEquitable}If $\pi=(d_{1},\ldots,d_{n})$ is a positive non-increasing degree sequence and some $G\in \mathcal{R}(\pi)$ has edge-disjoint regular factors $\{F_{1},\ldots,F_{p}\}$, then for \begin{equation}\label{eq:gammaBound}\gamma\geq \max_{l\leq m(\pi)}\bigg\{\bigg\lfloor\frac{d_{l}+l}{2}\bigg\rfloor\bigg\}+1\end{equation} there exists an $H\in \mathcal{R}(\pi)$ that
\begin{itemize}
    \item is equitably $\gamma$-colorable,
    \item has edge-disjoint regular factors $\{F'_{1},\ldots,F'_{p}\}$ such that $F'_{i}\in \mathcal{R}(F_{i})$, and
    \item if $\sum_{i=1}d_{i}\geq 2(n-1)$, then $H$ is $d_{n}-1$ edge-connected and $d_{n}$ edge-connected when $d_{n}$ is even or one.
\end{itemize}
\begin{proof}See Section~\ref{sec:equitablelColorable}.
\end{proof}
\end{theorem}
Since $d_{1}\geq d_{i}\geq m(\pi)$ for all $i\leq m(\pi)$, (\ref{eq:gammaBound}) is bounded above by $\bigg\lfloor\frac{d_{1}+m(\pi)}{2}\bigg\rfloor+1\leq d_{1}+1$ with equality when $d_{1}=m(\pi)$. On the other hand, equation (\ref{eq:gammaBound}) is sharp. To see this, consider the split graph $G=K^{s}*I^{2t-1}$ for $t\geq 1$. Note that $G$ is the unique realization of $\pi(G)$, the strong index of $\pi(G)$ is $s$, the first $s$ terms of $\pi(G)$ have degree $s + 2t - 2$, and $G$ is equitably $(s + t)$-colorable.

Separately, our two main theorems do not quite prove there is a realization that satisfies the Chen-Lih-Wu conjecture when $d_{1}=m(\pi)$ and no realization of $\pi$ has a regular factor. We will need the following theorem to connect the two.

\begin{theorem}\label{cor:bestBound}Let $\pi=(d_{1},\ldots,d_{n})$ be a non-increasing graphic sequence. For a non-negative integer $k\leq d_{n}$ such that $kn$ is even, if \begin{equation}\label{eq:bestBound}d_{d_{1}-d_{n}+1}\geq d_{1}-d_{n}+k-1,\end{equation} then some realization of $\pi$ has a $k$-factor.
\begin{proof}
    For even $n$,  Theorem~\ref{cor:bestBound} was first proved by Shook in \cite{Shook}. Using a different method we prove the theorem for all $n$ in Section~\ref{sec:kFactor}.
\end{proof}
\end{theorem}

 The bound in (\ref{eq:bestBound}) is sharp. For even $k$, let $G=I^{2}*H$ where $H$ is a $(k-2)$-regular graph with $k+2$ vertices. For odd $k$, let $G=K^{2}*H$ where $H$ is a $(k-2)$-regular graph with $k+1$ vertices. Thus, $\pi(G)=(d_{1},d_{2},\ldots, d_{k+4})$ is graphic where $d_{1}=d_{2}=k+2$ and $d_{3}=\ldots =d_{k+4}=k$, and $D_{k}(\pi(G))=(2,2,\ldots,0)$ is not graphic. Moreover, $d_{d_{1}-d_{k+4}+1}=d_{3}=k=d_{1}-d_{k+4}+k-2$.

We are now ready to prove a degree sequence version of the Chen-Lih-Wu conjecture.

\begin{theorem}\label{thm:chenLihWu}For $k\geq 0$, suppose $\pi=(d_{1},\ldots,d_{n})$ is a positive non-increasing degree sequence and some $G\in \mathcal{R}(\pi)$ has a $k$-factor that can be partitioned into edge-disjoint regular factors $\{F_{1},\ldots,F_{p}\}$. If $d_{n}\neq n-1$, $d_{1}\neq 1$, and $n$ is even when $d_{1}=d_{n}=2$, then there is an $H\in \mathcal{R}(\pi)$ that 
\begin{itemize}
    \item is equitably $d_{1}$-colorable,
    \item has edge-disjoint regular factors $\{F'_{1},\ldots,F'_{p}\}$ such that $F'_{i}\in \mathcal{R}(F_{i})$, and
    \item if $\sum_{i=1}d_{i}\geq 2(n-1)$, then $H$ is $d_{n}-1$ edge-connected and $d_{n}$ edge-connected when $d_{n}$ is even or one.
\end{itemize}
\begin{proof}If some realization of $\pi$ has a regular factor, then we may, without loss of generality, assume $F_{1}$ is not empty. If $m(\pi)<d_{1}$, then the theorem follows from Theorem~\ref{thm:ExistDegreeEquitable}. When $d_{1}=m(\pi)$, then since $d_{d_{1}-d_{n}+1}\geq d_{d_{1}}\geq d_{1}\geq d_{1}-d_{n}+k-1$, Theorem~\ref{cor:bestBound} says some realization of $\pi$ has a regular factor. Thus, $F_{1}$ is not empty, and the degree conditions on $\pi$ imply $G$ is not a complete graph, a $1$-factor, nor a $2$-factor when $n$ is odd. Therefore, the theorem follows from Theorem~\ref{thm:MaxDegreeEquitable} using $G$ and $F_{1}$. 
\end{proof}
\end{theorem}

\section{Our Past Work}\label{sec:PastWork}
Along with some new insights, the proofs in this paper rely on our results in \cite{Shook2024, Shook}. 

In Section~\ref{sec:ColorExchanges}, we present generalized edge-exchanges that we first studied in \cite{Shook}. These edge exchanges allow us to modify realizations so that the resulting graph still has a $k$-factor. We present a key lemma in that section that is critical to the proof of Theorem~\ref{thm:ExistDegreeEquitable}. 

In Section~\ref{sec:edgeconnectivity}, by simple modifications of our proofs in \cite{Shook2024} we show that we can modify an equitably $l$-colorable realization $H$ that has a $k$-factor $F$ so that the resulting realization $G$ is connected, equitably $l$-colorable, and has $H-E(F)$ as a subgraph.

\begin{theorem}\label{thm:connectednew}If there is an equitable $l$-colorable graph $G_{0}=(V,E)$ with spanning subgraph $Z_{0}$  such that $\delta(Z_{0})\geq 1$ and $|E(Z_{0})|\geq |V|-1$ when $\delta(G_{0})=1$, then there is an equitable $l$-colorable realization $G\in \mathcal{R}(G_{0},G_{0}-E(Z_{0}))$ such that $G$ is $\delta(G)-1$ edge-connected when $\delta(G)\geq 3$ and odd or $G$ is maximally edge-connected, otherwise.
\end{theorem}

Note that if one where to improve Theorem~\ref{thm:connectednew} so that $G$ is maximally edge-connected in all cases, then we can improve the main theorems of this paper.

\subsection{Generalized Edge-Exchanges}\label{sec:ColorExchanges}
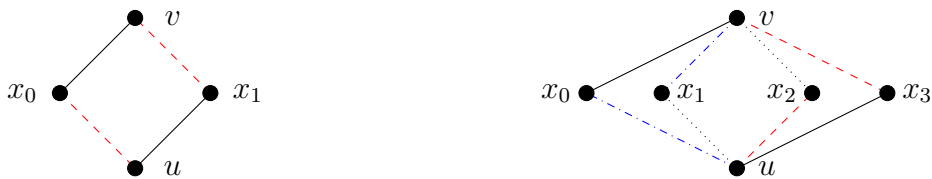
\begin{figure}[h]
    \centering
    \begin{tikzpicture}
\draw (1,3) -- (2,4);
\draw (3,3) -- (2,2);
\draw[dashed, color=red] (1,3) -- (2,2);
\draw[dashed, color=red] (3,3) -- (2,4);
\draw[fill=black] (1,3) circle (.1cm);
\draw[fill=black] (3,3) circle (.1cm);
\draw[fill=black] (2,2) circle (.1cm);
\draw[fill=black] (2,4) circle (.1cm);
\node[] at (.5,3) {$x_{0}$};
\node[] at (3.5,3) {$x_{1}$};
\node[] at (2.5,4) {$v$};
\node[] at (2.5,2) {$u$};


\draw (10,4) -- (8,3);
\draw (10,2) -- (12,3);
\draw[dashdotted, color=blue] (10,2) -- (8,3);
\draw[dashdotted, color=blue] (10,4) -- (9,3);
\draw[dotted, color=black] (10,2) -- (9,3);
\draw[dotted, color=black] (10,4) -- (11,3);
\draw[dashed, color=red] (10,4) -- (12,3);
\draw[dashed, color=red] (10,2) -- (11,3);

\draw[fill=black] (10,4) circle (.1cm);
\draw[fill=black] (10,2) circle (.1cm);
\draw[fill=black] (9,3) circle (.1cm);
\draw[fill=black] (8,3) circle (.1cm);
\draw[fill=black] (11,3) circle (.1cm);
\draw[fill=black] (12,3) circle (.1cm);

\node[] at (10.4,4) {$v$};
\node[] at (10.4,2) {$u$};
\node[] at (7.6,3) {$x_{0}$};
\node[] at (9.4,3) {$x_{1}$};
\node[] at (10.6,3) {$x_{2}$};
\node[] at (12.4,3) {$x_{3}$};

\end{tikzpicture}
    \caption{Edge-exchanges with length 2 and length 4, respectively. }
    \label{fig:Exchanges}
\end{figure}

We may transform one realization $G$ of a degree sequence $\pi$ to another by performing an operation called an edge-exchange or an edge-switch. This operation, that was first presented in 1891 by Petersen in \cite{Petersen1891}, removes two edges $xy$ and $uv$ of $G$ and replaces them with two non-edges $xu$ and $vy$ of $G$ (See the left side of Figure~\ref{fig:Exchanges}. This figure was reproduced from our paper \cite{Shook}). The resulting graph is still a realization of $\pi$. However, we need a more general operation that exchanges edges so that the presence of a regular factor is preserved. Seacrest in \cite{Seacrest2021} and Shook in \cite{Shook} explored this concept. For our purposes, we will use a more strict definition of edge-exchange than given in \cite{Shook}. 

Let $\pi=(d_{1},\ldots,d_{n})$ be a non-increasing degree sequence. Suppose $H$ is a realization of $\pi$ that has a $k$-factor that can be partitioned into $\{H_{3},\ldots, H_{p}\}$ regular factors. We let $H_{1}= H-\bigcup_{3\leq i\leq p} E(H_{i})$ and $H_{2}=\overline{H}$. To help with our discussion, we will assume for each $i$ that the edges of $H_{i}$ are colored $i$. We call a list \[L=(vx_{0},x_{0}u,vx_{1},x_{1}u,\ldots,vx_{q-1},x_{q-1}u)\] of at most $q\leq 2(p+2)$ edges such that for $j$ module $q$ the edges $x_{j}u$ and  $vx_{j+1}$ have the same color and $x_{s}u$ is a different color for all $s\neq j$ an edge-exchange. Indeed if we exchange the colors of $vx_{j}$ with $x_{j}u$ for all $j$ module $q$ we create a realization $H'$ of $\pi$ that has a $k$-factor that can be partitioned into $\{H'_{3},\ldots, H'_{p}\}$ regular factors such that $H'_{j}\in \mathcal{R}(H_{j})$ for $3\leq j \leq q$, and $H'_{1}\in \mathcal{R}(H_{1})$ where $H'_{1}= H'-\bigcup_{3\leq i\leq p} E(H'_{i})$. Let $\mathcal{X}(L)=\{x_{0},\ldots,x_{q-1}\}$. 

Shook \cite{Shook} gave conditions for when one can find internal vertex disjoint edge-exchanges. This ability plays a crucial role in the proof of Theorem~\ref{thm:ExistDegreeEquitable}. However, we don't need the full result. Thus, we end this section by presenting a simplified version of Lemma~7 in \cite{Shook} from the perspective of $\overline{H}$.

\begin{lemma}[\cite{Shook}]\label{lem:largedegree}Let $\pi=(d_{1},\ldots,d_{n})$ be a non-increasing degree sequence. Suppose $H$ is a realization of $\pi$ that has a $k$-factor that can be partitioned into $\{H_{3},\ldots, H_{p}\}$ regular factors. For vertices $u$ and $v$, let $X=\{x^{(1)}_{0},\ldots,x^{(|X|)}_{0}\}$ where $X\subseteq N_{\overline{H}}(v)-N_{\overline{H}}(u)$. If $deg_{H}(v)\geq deg_{H}(u)$ then there exists a set $\mathcal{L}=\{L^{(1)},\ldots,L^{(|X|)}\}$ of edge-exchanges such that for all $j$, $x^{(j)}_{0}\in \mathcal{X}(L^{(j)})$ and $\mathcal{X}(L^{(j)})\cap \mathcal{X}(L^{(i)})=\emptyset$ for $j\neq i$.
\end{lemma}

\subsection{Edge-Conectivity}\label{sec:edgeconnectivity}

Since the proof of Theorem~\ref{thm:connectednew} is similar to the proof found in \cite{Shook2024} of the following Theorem~\ref{thm:connectedLeftoverFull}, we are going to sketch the subtle differences and refer the reader to that paper. 

\begin{theorem}\label{thm:connectedLeftoverFull}[Theorem~2 in \cite{Shook2024}]If there is a graph $G_{0}=(V,E)$ with edge-disjoint spanning subgraphs $F$ and $Z_{0}$ with $\delta(Z_{0})>\Delta(F)$ such that $|E(Z_{0})|\geq |V|-1$ when $\delta(G_{0})=1$, then there is a $G\in \mathcal{R}(G_{0},G_{0}-E(Z_{0}))$ such that $G-E(F)$ is maximally edge-connected.
\end{theorem}

The parameter $F$ found in Theorem~\ref{thm:connectedLeftoverFull} is not needed in Theorem~\ref{thm:connectednew} so the reader can assume $F$ is empty when reading through. The reader will only need to follow the proof of Theorem~\ref{thm:connectedLeftoverFull} through Claim~8.3 in \cite{Shook2024}.  The proof begins by carefully choosing a counter-example $G$. Later two edges $aa'$ and $bb'$ of $Z_{0}$ are carefully chosen such that $E_{G}(\{a,a'\},\{b,b'\})=\emptyset$. Those edges are then exchanged with the two non-edges $ab'$ and $a'b$. However, for Theorem~\ref{thm:connectednew}, once those edges are exchanged, we need the resulting graph to be $l$ colorable. Claim~8.1 in \cite{Shook2024} shows that $G$ is connected and it is easy to follow. For the rest of the proof of Theorem~\ref{thm:connectedLeftoverFull}, two special sets $A$ and $B$ of $V(G)$ and two vertices $a\in A$ and a $b\in B$ such that $E_{G}(\{a\},\{b\})=\emptyset$ are under consideration. Moreover, in this situation they assume $E_{G}(A,B)<\delta(G)$. However, for Theorem~\ref{thm:connectednew} we replace this with the assumption that $E_{G}(A,B)<\lambda(Z_{0})$. This allows us to deduce that $a$ is adjacent in $Z_{0}$ to a vertex $a'\in A$ and that $b$ must be adjacent in $Z_{0}$ to a vertex $b'\in B$ such that $E_{G}(\{a,a'\},\{b,b'\})= \emptyset$. Now we can choose an arbitrary equitable $l$ coloring $f$ of $G$, and we may assume without loss of generality that $f(a)\neq f(b')$ and  $f(a')\neq f(b)$. This implies the resulting graph $G-\{aa',bb'\}+\{ab',ba'\}$ after the exchange has an equitable coloring $f'$ where $f(v)=f'(v)$ for all $v\in V(G)$. From here, the proof is nearly identical to that of Theorem~\ref{thm:connectedLeftoverFull}.

We suspect that Theorem~\ref{thm:connectednew} can be strengthened by proving that $G$ is maximally edge-connected in the odd case as well. The difficulty in proving this can be seen in Claim~8.6 in \cite{Shook2024}, where one of two possible edge exchanges must be performed, but we cannot guarantee the resulting realization will still be an equitable $l$ coloring for both edge-exchanges. This issue would have to be addressed.

\section{Proof of Theorem~\ref{thm:MaxDegreeEquitable}}\label{sec:chenlihWu}

Recall that a graph $H$ is said to pack with a graph $G$ if $H$ is a subgraph of the complement of $G$.

\begin{theorem}[Katerinis \cite{Katerinis1985}]\label{katerinisRegularCase}
    If $G_{1}$ and $G_{2}$ are simple graphs of order $n$ such that
    \begin{itemize}
        \item $G_{1}$ is $k$-regular for $k\geq 1$
        \item $n\geq 2(\Delta(G_{2})+1)$, and 
        \item $n\geq 4k-5$,
    \end{itemize}
then some realization of $\pi(G_{1})$ packs with $G_{2}$.
\end{theorem}

\begin{lemma}\label{lem:buildEquitable}For $1\leq k\leq \frac{n+5}{2}$ and a positive integer $i$, if $G=(V,E)$ is a simple $n$ vertex graph such that $G$ has a $k$-factor $F$ such that \[4-i\leq \Delta(G)-k\leq \frac{n}{4}-1+j\] where $j=1$ unless $\Delta(G)-k+i$ divides $n$, in which case $j=0$, then for $l\geq \Delta(G)-k+i$, there is a realization $H\in \mathcal{R}(\pi(G), G-E(F))$ that is equitably $l$-colorable. In addition, 
    we can require $H$ to be $(\delta(G)-1)$ edge-connected or $\delta(G)$ edge-connected when $\delta(G)$ is even or one. 
    \begin{proof}Let $Z=G-E(F)$ and $\Delta_{1}=\Delta(G)$.  Since $\Delta(Z)<\Delta_{1}$, Theorem~\ref{thm:HSTheorem} says there is an equitably $\Delta$-coloring $f$ of $Z$. Let $G_{2}$ be a graph on $V$ such that $G_{2}[Y_{i}(f)]$ is a clique for all $i$ and there are no other edges. Let $G'_{2}=Z+G_{2}$, and note that $\Delta(G_{2}')=\Delta_{1}-k+\bigg\lfloor\frac{n}{l}\bigg\rfloor-j$. By the assumptions in the statement of the lemma we may deduce that \[2(\Delta(G_{2}')+1)=2\left(\Delta_{1}-k+\bigg\lfloor\frac{n}{l}\bigg\rfloor-j+1\right)\leq  2\left(\frac{n}{4}-1+j +  \frac{n}{4}-j+1\right)= n.\] Therefore, by Theorem~\ref{katerinisRegularCase}, some realization of $F$ packs with $G_{2}'$.  We have now shown that there is a realization $H$ of $G$ that is equitably $l$-colorable and has $Z$ as a subgraph. However, in addition, if $H$ has at least $n-1$ edges and we let $F=Z_{0}$ in Theorem~\ref{thm:connectednew}, then we can require $H$ to be $(\delta(G)-1)$ edge-connected or $\delta(G)$ edge-connected when $\delta(G)$ is even or one.\qedhere
    \end{proof}
\end{lemma}

We will use the following two theorems to prove the main result of this section.

\begin{theorem}[\cite{Kierstead2012}]\label{thm:mD4}
    If $G$ is a $4$-colorable graph with $\Delta(G)\leq 4$, then $G$
has an equitable $\Delta(G)$-coloring.
\end{theorem}

\begin{theorem}[\cite{Kierstead2015}]\label{thm:Dgreaternover4}
    Let G be a graph with $\max\{\chi(G),\Delta(G),\frac{|G|}{4}\}\leq r$. If $r$ is even or $G$
does not contain $K^{r,r}$, then $G$ has an equitable r-coloring.
\end{theorem}

We restate Theorem~\ref{thm:MaxDegreeEquitable} for the reader's convenience.

\begin{reptheorem}{thm:MaxDegreeEquitable}For $k\geq 1$, if $G=(V,E)$ is a simple $n$ vertex graph such that $G$ has a $k$-factor $F$ and $G$ is not a complete graph, a $1$-factor, nor a $2$-factor when $n$ is odd, then there is a realization $H\in \mathcal{R}(\pi(G), G-E(F))$ that is equitably $\Delta(G)$-colorable. In addition, we can require $H$ to be $(\delta(G)-1)$ edge-connected or $\delta(G)$ edge-connected when $\delta(G)$ is even or one. 
\begin{proof}Let $Z=G-E(F)$, $\Delta_{1}=\Delta(G)$,  and note that $n=\Delta_{1}\bigg\lfloor\frac{n}{\Delta_{1}}\bigg\rfloor+r$ for some non-negative integer $r$.  We first show there is a realization $H\in \mathcal{R}(\pi(G), G-E(F))$ that is equitably $\Delta(G)$-colorable. If $4< \Delta_{1}< \bigg\lfloor\frac{n}{4}\bigg\rfloor$, then the theorem follows from an direct application of Lemma~\ref{lem:buildEquitable}. Suppose $\Delta_{1}\leq 4$. If $G$ is a $2$-factor with even $n$, then consecutively exchanging edges of $F$ in distinct odd cycles creates a realization of  $\mathcal{R}(\pi(G), G-E(F))$ with no odd cycles. Since there is some realization $H\in \mathcal{R}(\pi(G), G-E(F))$ that is not complete, a $1$-factor, nor a $2$-factor whith odd cycles, Theorem~\ref{thm:mD4} says $H$ is equitably $\Delta_{1}$-colorable.  We are left with the case $\Delta_{1}\geq \frac{n}{4}$. Suppose $\Delta_{1}$ is odd and $G=K^{\Delta_{1},\Delta_{1}}$. Let $A$ and $B$ be the two independent sets that makeup $G$. Let $xy$ and $uv$ be edges of $F$ such that $\{x,u\}\subseteq A$ and $\{v,y\}\subseteq B$. If we remove the edges $xy$ and $uv$ from $G$ and add the edges $xu$ and $yv$ to create a new graph $G'$, then $G'$ contains $Z$ and is no longer a complete bipartite graph since $G$ is not a four-cycle. Thus, there is a realization $H$ of $\pi(G)$ that contains $Z$ and is not $K^{\Delta_{1},\Delta_{1}}$ when $\Delta_{1}$ is odd. Theorem~\ref{thm:Dgreaternover4} says that $H$ can be equitably $\Delta_{1}$-colorable. We have now shown that there is a realization $H$ of $G$ that is equitably $\Delta_{1}$-colorable and has $Z$ as a subgraph. However, in addition, if $H$ has at least $n-1$ edges and we let $F=Z_{0}$ in Theorem~\ref{thm:connectednew}, then we can require $H$ to be $(\delta(G)-1)$ edge-connected or $\delta(G)$ edge-connected when $\delta(G)$ is even or one.\qedhere
    
\end{proof}
\end{reptheorem}

\section{Proof Of Theorem~\ref{thm:ExistDegreeEquitable}}\label{sec:equitablelColorable}

\begin{reptheorem}{thm:ExistDegreeEquitable}If $\pi=(d_{1},\ldots,d_{n})$ is a positive non-increasing degree sequence and some $G\in \mathcal{R}(\pi)$ has edge-disjoint regular factors $\{F_{1},\ldots,F_{p}\}$, then for \begin{equation}\gamma\geq \max_{l\leq m(\pi)}\bigg\{\bigg\lfloor\frac{d_{l}+l}{2}\bigg\rfloor\bigg\}+1\end{equation} there exists an $H\in \mathcal{R}(\pi)$ that
    \begin{enumerate}[label=(E\arabic*),ref=(E\arabic*)]
    \item \label{main:req1}is equitably $\gamma$-colorable,
        \item \label{main:req2}has edge-disjoint regular factors $\{F'_{1},\ldots,F'_{p}\}$ such that $F'_{i}\in \mathcal{R}(F_{i})$, and
        \item \label{main:req3} if $\sum_{i=1}d_{i}\geq 2(n-1)$, then $H$ is $d_{n}-1$ edge-connected and $d_{n}$ edge-connected when $d_{n}$ is even or one.
    \end{enumerate}
    
    \begin{proof}Every realization of $\pi(G)$ is equitably $n$ colorable. Therefore, by contradiction, we assume there is a largest $\gamma$ that does not satisfy the theorem. Thus, we can carefully select a realization with the required properties that can be $\gamma+1$ equitably colored. We let $m:=m(\pi)$. Let $\mathcal{H}\subseteq \mathcal{R}(\pi)$ be the largest set where every $H\in \mathcal{H}$ has a set of edge-disjoint regular factors $\{F'_{1},\ldots,F'_{p}\}$ such that $F'_{i}\in \mathcal{R}(F_{i})$. We choose an $H\in \mathcal{H}$ such that
    \begin{enumerate}[label=(\roman*),ref=(\roman*)]
        \item \label{main:choice1}there is an $\gamma+1$ coloring $f$ of $H$ where $|Y_{1}(f)|\geq \ldots \geq |Y_{\gamma+1}(f)|$ and  $|Y_{1}(f)|\geq |Y_{i}(f)|\geq |Y_{1}(f)|-1$ for $i\leq \gamma$,
        \item \label{main:choice2} subject to \ref{main:choice1}, we  minimize $|Y_{\gamma+1}(f)|$, and
        \item \label{main:choice3} subject to \ref{main:choice2}, for $2\leq j\leq \gamma+1$, we minimize in order the sum \[\sum_{i=1}^{j-1}\sum_{v_{z}\in Y_{i}(f)}z.\] 
        \item \label{main:choice4} subject to \ref{main:choice3}, we minimize \[\sum_{i=1}^{\gamma}\sum_{v_{z}\in N_{H}(v_{\alpha_{\gamma+1}})\cap Y_{i}(f)}z\] where $\alpha_{\gamma+1}=\max\{z|v_{z}\in Y_{\gamma+1}(f)\}$.
        
    \end{enumerate}
        We may assume $Y_{\gamma+1}(f)\neq \emptyset$ since otherwise, $f$ would be an equitable $\gamma$ coloring of $H$. Let $q$ be the smallest index such that $|Y_{q}(f)|=|Y_{\gamma}(f)|$.
    
    \begin{claim}\label{cl:basicNeighborhood1}Every vertex in $Y_{\gamma+1}(f)$ is adjacent in $H$ to a vertex in $Y_{i}(f)$ for $\gamma\geq i\geq q$.
    
    \begin{proof}If there were a vertex $x\in Y_{\gamma+1}(f)$ not adjacent in $H$ to a vertex in $Y_{i}(f)$ for $l\geq i\geq q$, then we may remove $x$ from $Y_{\gamma+1}(f)$ and add $x$ to $Y_{i}(f)$ to create a new coloring of $H$ that violates \ref{main:choice2}. 
    \end{proof}
    \end{claim}
    
    \begin{claim}\label{cl:basicNeighborhood2}For a $v_{i}\in Y_{\gamma+1}(f)$ and a $Y_{j}(f)$ with $\gamma\geq j\geq q$, every vertex $v_{s}\notin Y_{j}(f)$ with $d_{s}\geq d_{i}$ is adjacent in $H$ to some vertex in $Y_{j}(f)$.
    \begin{proof}Let $N_{H}(v_{i})\cap Y_{j}(f)=\{x_{1},x_{2},\ldots,x_{p}\}$, and suppose some vertex $v_{s}\in Y_{t}(f)$ for some $t\neq j$ is not adjacent in $H$ to a vertex in $Y_{j}(f)$. Since $d_{s}\geq d_{i}$ and $v_{s}$ is not adjacent in $H$ to vertices in $\{x_{1},\ldots,x_{p}\}$, Lemma~\ref{lem:largedegree} says we can find $p$ edge-exchanges $\mathcal{L}=\{L^{(1)},\ldots,L^{(p)}\}$ between $v_{i}$ and $v_{s}$ such that $x_{z}\in \mathcal{X}(L^{(z)})$ and $\mathcal{X}(L^{(z)})\cap \mathcal{X}(L^{(z')})=\emptyset$ for $z\neq z'$. Since no edge-exchange in $\mathcal{L}$ share internal vertices, we can perform all $p$ edge-exchanges to construct a realization $H'$ of $\pi(G)$. Note that $v_{i}$ is not adjacent in $H'$ to vertices in $Y_{j}(f)$ and $v_{s}$ is not adjacent in $H'$ to the vertices in $Y_{t}(f)-\{v_{s}\}$. Thus, there is a $\gamma+1$ coloring $f'$ of $H'$ with color classes $Y_{z}(f')=Y_{z}(f)$ for $z\notin\{j,\gamma+1\}$, $Y_{j}(f')=Y_{j}(f)+\{v_{i}\}$, and $Y_{\gamma+1}(f')=Y_{\gamma+1}(f)-\{v_{i}\}$. However, this violates \ref{main:choice2}.
    \end{proof}
    \end{claim}
    
    \begin{claim}\label{cl:basicDegreeBound}For $1\leq z\leq z'-1< \gamma+1$ and $v_{t}\in Y_{z}(f)$, if some $v_{s}\in Y_{z'}(f)$ is not adjacent in $H$ to vertices in $Y_{z}(f)-\{v_{t}\}$, then $t<s$.
    \begin{proof}By contradiction we assume $t>s$. Let $N_{H}(v_{t})\cap (Y_{z'}(f)-\{v_{s}\})=\{x_{1},x_{2},\ldots,x_{p}\}$. Since $d_{s}\geq d_{t}$ and $v_{s}$ is not adjacent in $H$ to vertices in $\{x_{1},\ldots,x_{p}\}$, Lemma~\ref{lem:largedegree} says we can find $p$ edge-exchanges $\mathcal{L}=\{L^{(1)},\ldots,L^{(p)}\}$ between $v_{t}$ and $v_{s}$ such that $x_{i}\in \mathcal{X}(L^{(i)})$ and $\mathcal{X}(L^{(i)})\cap \mathcal{X}(L^{(j)})=\emptyset$ for $i\neq j$. Moreover, $\mathcal{X}(L^{(i)})\cap (Y_{z}(f)-v_{t})=\emptyset$ since $v_{s}$ is not adjacent to the vertices in $\{x_{1},x_{2},\ldots,x_{p}\}$. Therefore, we may perform all of the $p$ edge-exchanges to construct another realization $H'$ of $\pi(G)$. Note that $v_{t}$ is not adjacent in $H'$ to vertices in $Y_{z'}(f)-\{v_{s}\}$ and $v_{s}$ is not adjacent in $H'$ to the vertices in $Y_{z}(f)-\{v_{t}\}$. Thus, there is a $\gamma+1$ coloring $f'$ of $H'$ with color classes $Y_{i}(f')=Y_{i}(f)$ for $i\notin\{z, z'\}$, $Y_{z}(f')=Y_{z}(f)-\{v_{t}\}+\{v_{s}\}$, and $Y_{z'}(f')=Y_{\gamma+1}(f)-\{v_{s}\}+\{v_{t}\}$. However, this violates \ref{main:choice4} since  \[\sum_{i=1}^{z'-1}\sum_{v_{j}\in Y_{i}(f')}j<\sum_{i=1}^{z'-1}\sum_{v_{j}\in Y_{i}(f)}j.\qedhere\]
    \end{proof}
    \end{claim}
    
    Let $\beta_{\gamma}=\min\{z| v_{z}\in N_{H}(v_{\alpha_{\gamma+1}})\cap Y_{\gamma}(f)\}$, and $\alpha_{i}=\max\{j|v_{j}\in Y_{i}(f)\}$ for all $Y_{i}$.  Claim~\ref{cl:basicDegreeBound} implies that if $v_{\alpha_{i}}$ is not adjacent to a vertex in $Y_{j}(f)$ for $j\leq i-1$, then $z<\alpha_{i}$ for every $v_{z}\in Y_{j}(f)$.
    
    \begin{claim}\label{cl:lowerDegreebound}If $d_{\alpha_{\gamma+1}}<\gamma$, then $d_{\beta_{\gamma}}\geq \gamma$.
        \begin{proof}Let $N_{H}(v_{\alpha_{\gamma+1}})\cap Y_{\gamma}=\{x_{1},x_{2},\ldots,x_{p}\}$. We partition the set $Y_{1}(f),\ldots, Y_{\gamma-1}(f)$ into two sets $A$ and $B$ such that every $Y_{i}(f)$ that has a $v_{z}\in Y_{i}$ with $z>\alpha_{\gamma+1}$ is in $B$ and the rest are in $A$. Let $\hat{A}=\bigcup_{Y_{z}(f)\in A}Y_{z}$. By Claim~\ref{cl:basicDegreeBound} $v_{\alpha_{\gamma+1}}$ is adjacent to at least one vertex in each set of $B$. If $v_{\alpha_{\gamma+1}}$ is adjacent to at least $|A|-(p-1)$ sets in $A$, then $d_{\alpha_{\gamma+1}}\geq |A|-(p-1)+|B|+|\{x_{1},x_{2},\ldots,x_{p}\}|=\gamma-1-(p-1)+p\geq \gamma$. We now assume $v_{\alpha_{\gamma+1}}$ is not adjacent to vertices in at least $p$ sets in $A$. By Claim~\ref{cl:basicNeighborhood1}, each of those sets must have $|Y_{\gamma}(f)|+1$ vertices. 
    
        By Claim~\ref{cl:basicNeighborhood2} every vertex in $\hat{A}$ is adjacent in $H$ to a vertex in $Y_{\gamma}(f)$. This implies \[e_{H}(\hat{A},Y_{\gamma}(f))\geq |\hat{A}|\geq |A||Y_{\gamma}(f)|+p.\] Thus, some vertex $v_{s}\in Y_{\gamma}(f)$ must be adjacent in $H$ to at least $|A|+1$ vertices in $\hat{A}$. First, we consider the case $s<\alpha_{\gamma+1}$. By Claim~\ref{cl:basicDegreeBound}, $v_{s}$ is adjacent in $H$ to at least one vertex in each set of $B$. This implies $d_{s}\geq |A|+1+|B|=\gamma$. We are left with two subcases. If $\beta_{\gamma}<s$, then $d_{\beta_{\gamma}}\geq \gamma$ and we are done. If $\beta_{\gamma}>s$, then by our choice of $v_{\beta_{\gamma}}$, $v_{s}$ is not adjacent in $H$ to $v_{\alpha_{\gamma+1}}$. However, we may exchange the edge $v_{\beta_{\gamma}}v_{\alpha_{\gamma+1}}$ with the non-edge $v_{\alpha_{\gamma+1}}v_{s}$ to construct a realization that contradicts \ref{main:choice4}. So we are left with the case that $s>\alpha_{\gamma+1}$. Let $(N_{H}(v_{s})\cap \hat{A})-N_{H}(v_{\alpha_{\gamma+1}})=\{x'_{1},\ldots, x'_{\lambda}\}$. Lemma~\ref{lem:largedegree} says we can find $\lambda$ edge-exchanges $\mathcal{L}=\{L^{(1)},\ldots,L^{(\lambda)}\}$ between $v_{s}$ and $v_{\alpha_{\gamma+1}}$ such that $x'_{j}\in \mathcal{X}(L^{(j)})$ and $\mathcal{X}(L^{(j)})\cap \mathcal{X}(L^{(i)})=\emptyset$ for $j\neq i$. Without loss of generality we may assume the first $\lambda'$ edge-exchanges $\mathcal{L'}=\{L^{(1)},\ldots, L^{(\lambda')}\}$ are such that $\mathcal{X}(L^{(i)})\cap Y_{\gamma}(f)=\emptyset$ for $i\leq \lambda'$. Since every edge-exchange in $\mathcal{L'}$ is internally vertex disjoint, we can perform all $\lambda'$ edge-exchanges to construct a realization $H'$ of $\pi(G)$. Thus, $v_{\alpha_{\gamma+1}}$ is still adjacent in $H'$ to the same $p$ neighbors in $Y_{\gamma}(f)$. By Claim~\ref{cl:basicDegreeBound} $v_{\alpha_{\gamma+1}}$ must be adjacent in $H'$ to every set in $B$. Note that $\lambda'\geq \lambda-p$. Thus, $v_{\alpha_{\gamma+1}}$ is adjacent in $H'$ to at least \[|N_{H}(v_{\alpha_{\gamma+1}})|+\lambda'=|(N_{H}(v_{s})\cap \hat{A})|-\lambda+\lambda'\geq |A|+1-\lambda+\lambda'\geq |A|+1-p\]  vertices in $\hat{A}$. Thus, $v_{\alpha_{\gamma+1}}$ is adjacent in $H'$ to at least $|A|+1-p+|B|+p=|A|+|B|+1=\gamma$ vertices in $H'$. Thus, $d_{\alpha_{\gamma+1}}\geq \gamma$. This completes the proof of the claim.
        \end{proof}
    \end{claim}
    
    By Claim~\ref{cl:lowerDegreebound} there is a smallest $s\in \{\alpha_{\gamma+1},\beta_{\gamma}\}$ such that $d_{s}\geq \gamma$. We now partition the sets $Y_{1}(f),\ldots, Y_{\gamma-1}(f)$ into three parts $P_{0}$, $P_{1}$, and $P_{2}$ where $v_{s}$ is not adjacent to any vertex in the sets of $P_{0}$, $v_{s}$ is adjacent to exactly one vertex in each of the sets of $P_{1}$, and is adjacent in $G$ to at least two vertices in each of the sets of $P_{2}$. By Claim~\ref{cl:basicDegreeBound} every vertex in the sets of $P_{0}$ must have a lower index than $s$. Therefore, if $|Y_{\gamma}(f)|=1$, then $P_{0}$ is empty. This is because Claim~\ref{cl:basicDegreeBound} says $s=\beta_{\gamma}<\alpha_{\gamma+1}$, and therefore, Claim~\ref{cl:basicNeighborhood2} says $v_{s}$ must be adjacent to every vertex that makes up the sets of $P_{0}$. Thus, $|Y_{\gamma}(f)|\geq 2$ when $P_{0}$ is not empty. Moreover, Claim~\ref{cl:basicDegreeBound} implies the sole neighbors of $v_{s}$ in each set of $P_{1}$ must have a lower index. Thus, \begin{equation*}
        s>|P_{0}||Y_{\gamma}(f)|+|P_{1}|\geq 2|P_{0}|+|P_{1}|.
    \end{equation*}
    Therefore, \begin{align}
        d_{s}&\geq |N_{H}(v_{s})\cap (Y_{\gamma}(f)\cup Y_{\gamma+1})(f)|+ |P_{1}|+2(\gamma+1-|P_{0}|-|P_{1}|-2)\notag\\
        &=|N_{H}(v_{s})\cap (Y_{\gamma}(f)\cup Y_{\gamma+1}(f))|+ 2\gamma-(2|P_{0}|+|P_{1}|)-2\notag\\
        &\geq |N_{H}(v_{s})\cap (Y_{\gamma}(f)\cup Y_{\gamma+1}(f))|+ 2\gamma-(s-1)-2\notag\\
        &= |N_{H}(v_{s})\cap (Y_{\gamma}(f)\cup Y_{\gamma+1}(f))|+ 2\gamma-s-1.\notag
    \end{align}
    Since $s\in \{\alpha_{\gamma+1},\beta_{\gamma}\}$ and $v_{\alpha_{\gamma+1}}v_{\beta_{\gamma}}$ is an edge of $H$, $|N_{H}(v_{s})\cap (Y_{\gamma}(f)\cup Y_{\gamma+1}(f))|\geq 1$. Thus, $\frac{d_{s}+s}{2}\geq \gamma$. If $s\leq m$, then \[\max_{l\leq m}\bigg\{\bigg\lfloor\frac{d_{l}+l}{2}\bigg\rfloor\bigg\}\geq \bigg\lfloor\frac{d_{s}+s}{2}\bigg\rfloor\geq \gamma.\]  If $s>m$, then \[\max_{l\leq m}\bigg\{\bigg\lfloor\frac{d_{l}+l}{2}\bigg\rfloor\bigg\}\geq \bigg\lfloor\frac{d_{m}+m}{2}\bigg\rfloor\geq m\geq d_{s}\geq \gamma.\] Thus, we have shown that $\gamma\leq \max_{l\leq m}\bigg\{\bigg\lfloor\frac{d_{l}+l}{2}\bigg\rfloor\bigg\}$. If $\sum_{i=1}^{n}d_{i}=2(n-1)$, then any realization $H\in \mathcal{R}(\pi)$ that satisfies \ref{main:req1} and \ref{main:req2} can be modified to satisfy \ref{main:req3}. To see this we simply apply Theorem~\ref{thm:connectednew} to $H$ by letting $Z_{0}$ be the $k$-factor formed by $\{F'_{1},\ldots,F'_{p}\}$ in $H$. Thus, proving our theorem.\qedhere
      
    \end{proof}
    \end{reptheorem}
        
\section{Proof of Theorem~\ref{cor:bestBound}}\label{sec:kFactor}

Our proof incorporates two theorems that are fundamental to the study of degree sequences. We start with an improvement to the seminal Erd\'{o}s-Gallai Theorem \cite{Erdos1960} that incorporates the strong index. 
\begin{theorem}[\cite{Zverovich1992,Li1975,Hammer1981}]\label{thm:ErdosGallai} A non-negative non-increasing sequence $\pi=(d_{1},\ldots,d_{n})$ is graphic if $\sum_{i=1}^{n}d_{i}$ is even and for all $l\leq m(\pi)$, \[\sum_{i=1}^{l}d_{i}\leq l(l-1)+\sum_{i=l+1}^{n}\min\{l,d_{i}\}.\] 
\end{theorem}
Note that Zverovich and Zverovich proved Theorem~\ref{thm:ErdosGallai} as stated with a simple observation that the theorem of Hammer, Ibaraki, Simeone \cite{Hammer1981}, and Li\cite{Li1975} only needed to check one less inequality.

Next, we enlist the help of Kundu's $k$-factor theorem.

\begin{theorem}[Kundu's $k$-factor Theorem \cite{Kundu1974} (See \cite{Chen1988} for a short proof.)]\label{kundu}
For $k\geq 0$, if $\pi=(d_{1},\ldots,d_{n})$ and $(d_{1}-k_{1},\ldots,d_{n}-k_{n})$ are both graphic such that $k\leq k_{i}\leq k+1$ for $1\leq i\leq n$, then there exists a realization of $\pi$ that has a $(k_{1},\dots,k_{n})$-factor.
\end{theorem}

For some non-negative integer $k$ and $\pi=(d_{1},\ldots, d_{n})$, we denote by $\mathcal{D}_{k}(\pi)$ the sequence $(d_{1}-k,\ldots,d_{n}-k)$. The following Lemma of A.R. Rao and S.B. Rao allows us to consider the largest $k'\geq k$ such that $k'n$ is even and $\mathcal{D}_{k'}(\pi)$ is not graphic.

\begin{lemma}[\cite{RamachandraRao1972}]\label{lem:ErdosRao}For non-negative integer $k$, let $\pi$ be a graphic degree sequence such that $D_{k}(\pi)$ is also graphic. If
$r$ is a positive integer such that $r \leq  k$ and $rn$ is even, then $D_{r}(\pi)$ is also graphic.
\end{lemma}

As a convenience to the reader, we restate Theorem~\ref{cor:bestBound}.

\begin{reptheorem}{cor:bestBound}Let $\pi=(d_{1},\ldots,d_{n})$ be a non-increasing graphic sequence. For a non-negative integer $k\leq d_{n}$ such that $kn$ is even, if \begin{equation*}d_{d_{1}-d_{n}+1}\geq d_{1}-d_{n}+k-1,\end{equation*} then some realization of $\pi$ has a $k$-factor.
\begin{proof}
We choose the largest $k'\geq k$ such that $k'n$ is even and $d_{d_{1}-d_{n}+1}\geq d_{1}-d_{n}+k'-1$. Thus, $d_{1}-d_{n}+k'\geq d_{d_{1}-d_{n}+1}\geq d_{1}-d_{n}+k'-1$. If $D_{k'}(\pi)=(d'_{1},\ldots, d'_{n})$ is graphic, then Kundu's Theorem implies some realization of $\pi$ has a $k'$-factor. Consequently, Lemma~\ref{lem:ErdosRao} implies some realization of $\pi$ has a $k$-factor. Thus, we assume $D_{k'}(\pi)$ is not graphic, and therefore, $d_{1}>d_{n}$.

We let $m=m(D_{k'}(\pi))$. Since $d_{1}-d_{n}\geq d_{d_{1}-d_{n}+1}-k'=d'_{d_{1}-d_{n}+1}$, $m\leq d_{1}-d_{n}$, and therefore, $n\geq\Delta_{1}+1\geq d_{1}-d_{n}+k+1\geq m+k+1$. 

Since $\sum_{i=1}^{n}d'_{i}=\sum_{i=1}^{n}d_{i}-k'n$ and both $\sum_{i=1}^{n}d_{i}$ and $k'n$ are even, $\sum_{i=1}^{n}d'_{i}$ is even. Thus, $D_{k'}(\pi)$ satisfies the first requirement of Theorem~\ref{thm:ErdosGallai}. 

Since $D_{k'}(\pi)$ is not graphic, Theorem~\ref{thm:ErdosGallai} says there is an $l\leq m$ such that \begin{equation}\label{eq:lowl}
    \sum_{i=1}^{l}d'_{i}>l(l-1)+\sum_{i=l+1}^{n}\min\{l,d'_{i}\}.
\end{equation} 
If $d'_{n-k'-1}\geq l$, then $\sum_{i=l+1}^{n}\min\{l,d'_{i}\}\geq (n-l-k')l$. This implies the contradiction \[\sum_{i=1}^{l}d'_{i}>l(l-1)+\sum_{i=l+1}^{n}\min\{l,d'_{i}\}\geq (n-1-k')l\geq  \sum_{i=1}^{l}d'_{i}.\] Thus, there is a smallest $t \leq n-(k'+1)$ such that $d'_{t}<l$. 

By definition $t\geq l+1$, and \[\sum_{i=t}^{n}\min\{l,d'_{i}\}=\sum_{i=t}^{n}d'_{i}=\bigg(\sum_{i=t}^{n}d_{i}\bigg)-(n-(t-1))k'=\bigg(\sum_{i=t}^{n}\min\{l,d_{i}\}\bigg)-(n-(t-1))k'.\] Moreover, if $t\geq l+2$, then $\sum_{i=l+1}^{t-1}\min\{l,d'_{i}\}=\sum_{i=l+1}^{t-1}\min\{l,d_{i}\}$. Since $\pi$ is graphic, \[\sum_{i=1}^{l}d_{i}\leq l(l-1)+\sum_{i=l+1}^{n}\min\{l,d_{i}\}.\] Therefore, from (\ref{eq:lowl}) we see that \[\sum_{i=1}^{l}d'_{i}>l(l-1)+\bigg(\sum_{i=l+1}^{n}\min\{l,d_{i}\}\bigg)-(n-(t-1))k'\geq \bigg(\sum_{i=1}^{l}d_{i}\bigg)-(n-(t-1))k'.\] Since $\sum_{i=1}^{l}d'_{i}=\bigg(\sum_{i=1}^{l}d_{i}\bigg)-lk'$, we may conclude that $(n-(t-1))k'>lk'$. Therefore, $n-t=l+\epsilon$ for some non-negative integer $\epsilon$.

We now consider the extremes of (\ref{eq:lowl}). To aid 
 us we let $c_{0}=d'_{t}+\bigg(\sum_{i=t+1}^{n}d'_{i}\bigg)-(l+\epsilon)(d_{n}-k')$, and $c_{1}=l(d_{1}-k')-\sum_{i=1}^{l}d'_{i}$. Note that $c_{0}\geq d'_{t}$ and $c_{1}\geq 0$, since $d_{1}-k'\geq d'_{i}$ and $n-t=l+\epsilon$.
We deduce that
\begin{align}
    l(d_{1}-k')-c_{1}=\sum_{i=1}^{l}d'_{i}&>l(l-1)+l(t-1-l)+\sum_{i=t}^{n}d'_{i}\notag\\
    &\geq l(t-2)+(l+\epsilon)(d_{n}-k')+c_{0}.\notag
\end{align}
Thus, \[l(d_{1}-k')-l(t-2)-c_{1}= l(d_{1}-k'-(t-2))-c_{1}>(l+\epsilon)(d_{n}-k')+c_{0}.\]
Therefore, 

\[c_{0}+c_{1}<l(d_{1}-k'+2-t-d_{n}+k')=l(d_{1}-d_{n}+2-t)-\epsilon(d_{n}-k').\] 
This implies $t\leq d_{1}-d_{n}+1$ since $c_{0}+c_{1}\geq 0$. Since $d_{1}-d_{n}\geq l>d'_{t}\geq d'_{d_{1}-d_{n}+1}\geq d_{1}-d_{n}-1$, $t=d_{1}-d_{n}+1$, $d'_{t}=d_{1}-d_{n}-1$, and $l=d_{1}-d_{n}=t-1$. Therefore, \[(d_{1}-d_{n})-\epsilon(d_{n}-k')>c_{0}+c_{1}\geq d'_{t}+c_{1}= d_{1}-d_{n}-1 + c_{1}.\] This implies $c_{0}=d'_{t}$ and $c_{1}=\epsilon(d_{n}-k')=0$. As a result, we may deduce that $\sum_{i=t+1}^{n}d'_{i}=(d_{1}-d_{n})(d_{n}-k')$ and $\sum_{i=1}^{t-1}d'_{i}=\sum_{i=1}^{l}d'_{i}=(d_{1}-d_{n})(d_{1}-k')$. This implies \[\sum_{i=1}^{n}d'_{i}=\sum_{i=1}^{t-1}d'_{i}+d'_{t}+\sum_{i=t+1}^{n}d'_{i}=(d_{1}-d_{n})(d_{1}-k')+d_{1}-d_{n}-1+(d_{1}-d_{n})(d_{n}-k').\] Collecting terms and rearranging we see that \[\sum_{i=1}^{n}d'_{i}=(d_{1}-d_{n})(d_{1}+d_{n}-2k'+1)-1=(d_{1}-d_{n})(d_{1}-d_{n}+1)+2(d_{1}-d_{n})(d_{n}-k')-1.\]

Since $(d_{1}-d_{n})(d_{1}-d_{n}+1)-1$ is odd, we have a contradiction to $\sum_{i=1}^{n}d'_{i}$ being even. Thus, $D_{k'}(\pi)$ is graphic.
\end{proof}  
\end{reptheorem}

\providecommand{\bysame}{\leavevmode\hbox to3em{\hrulefill}\thinspace}
\providecommand{\MR}{\relax\ifhmode\unskip\space\fi MR }
\providecommand{\MRhref}[2]{%
  \href{http://www.ams.org/mathscinet-getitem?mr=#1}{#2}
}
\providecommand{\href}[2]{#2}

\end{document}